\RequirePackage{fix-cm}
\documentclass[smallextended]{svjour3}       
\smartqed  
\usepackage{graphicx}
\usepackage{amsfonts}
\usepackage{color}
\usepackage{amsmath}
\definecolor{c20}{rgb}{0.,0.7,0.}
\definecolor{c30}{rgb}{0.,0.,1.}
\definecolor{c40}{rgb}{1,0.1,0.7}
\definecolor{c50}{rgb}{1,0,0}
\definecolor{c60}{rgb}{1,0.9,0.1}

\def\Ae#1{\textcolor{c20}{#1}}
\def\Ae#1{#1}

\def\eH#1{\textcolor{c30}{#1}}
\def\eH#1{#1}

\def\hH#1{\textcolor{c30}{#1}}
\def\hH#1{#1}

\def\cL#1{\textcolor{c50}{#1}}
\def\cL#1{#1}

\def\cling#1{\textcolor{c20}{#1}}
\def\cling#1{#1}

\newcommand{\tbb}[1]{{\textcolor{blue}{#1}}}
\def\tbb#1{#1}

\newcommand{\kb}[1]{\boldsymbol{#1}}
\newcommand{\vk}[1]{\kb{#1}}
\def\kal#1{{\cal{ #1}}}

\newcommand{\ve}{\varepsilon}
\def\fracl#1#2{\biggr(\frac{#1}{#2} \biggl) }

\newcommand{\abs}[1]{\left\lvert #1 \right\rvert}
\newcommand{\Abs}[1]{ \biggl \lvert #1 \biggr \rvert}
\newcommand{\ABs}[1]{ \biggl \lvert #1 \biggr \rvert}

\newcommand{\E}[1]{\mathbb{E}\left(#1\right)}

\newcommand{\pk}[1]{\mathbb{P} \left( #1 \right ) }

\newcommand{\R}{\mathbb{R}}

\newcommand{\N}{\mathbb{N}}
\newcommand{\inr}{\in \R}
\newcommand{\inn}{\in \N}
\newcommand{\ldot}{,\ldots,}

\newcommand{\limit}[1]{\lim_{#1 \to   \infty}}

\newcommand{\BQN}{\begin{eqnarray}}
\newcommand{\EQN}{\end{eqnarray}}
\newcommand{\BQNY}{\begin{eqnarray*}}
\newcommand{\EQNY}{\end{eqnarray*}}

\newcommand{\BS}{\begin{sat}}
\newcommand{\ES}{\end{sat}}

\newcommand{\BD}{\begin{de}}
\newcommand{\ED}{\end{de}}

\newcommand{\BIT}{\begin{itemize}}
\newcommand{\EIT}{\end{itemize}}
\newcommand{\BDI}{\begin{description}}
\newcommand{\EDI}{\end{description}}

\newtheorem{theo}{Theorem}[section]
\newtheorem{sat}[theo]{Proposition}
\newtheorem{de}[theo]{Definition}

\newcommand{\nelem}[1]{{Lemma \ref{#1}}}

\newcommand{\netheo}[1]{{Theorem \ref{#1}}}

\newcommand{\nesec}[1]{{Section~\ref{#1}}}

\newcommand{\COM}[1]{}

\def\cE#1{#1}

\def\IF{\infty}

\newcommand{\expon}[1]{\exp\left(#1\right)}

\def\cAL{\cE{\kal{C}_{\AL,\begin{lemma}}}}

\date{}

\def\X{\vk X}

\def\equivdis{\stackrel{d}{=}}
\def\ARIJ{A_{ij}}

\def\Cov{\mathbf{Cov}}
\newcommand{\todis}{\stackrel{d}{\to}}
\newcommand{\toprob}{\stackrel{p}{\to}}

\newcommand{\tb}[1]{{\textcolor{blue}{#1}}}

\def\tb#1{#1}
\newcommand{\tn}[1]{{\textcolor{blue}{#1}}}
\def\tn#1{#1}
\def\cL#1{#1}

\begin{document}
\title{Extremes of Order Statistics of Stationary Processes
\thanks{Partial support from the Swiss National Science Foundation Project 200021-140633/1  and Marie Curie International Research Staff Exchange Scheme Fellowship within the 7th European Community Framework Programme (Grant No. RARE-318984)
 is kindly acknowledged. The first author also acknowledges
 partial support by \tbb{Narodowe Centrum Nauki} Grant No 2013/09/B/ST1/01778 (2014-2016).}
}

\titlerunning{Extremes of order statistics of stationary processes}        

\author{Krzysztof D\c{e}bicki        \and
        Enkelejd Hashorva  \and
        Lanpeng Ji \and
        Chengxiu Ling
}
\institute{K. D\c{e}bicki \at
              Mathematical Institute, University of Wroc\l aw, \\ \cling{pl.} Grunwaldzki 2/4,
              50-384 Wroc\l aw, Poland
           \and
           E. Hashorva \and  L. Ji \and C. Ling\at
           Department of Actuarial Science,
University of Lausanne,\\
UNIL-Dorigny, 1015 Lausanne, Switzerland\\
 \email{chengxiu.ling@unil.ch}
}
\date{Received: date / Accepted: date}

\maketitle

\begin{abstract}
Let $\{X_i(t),t\ge0\}, 1\le i\le n$ be independent copies of a stationary  process $\{X(t), t\ge0\}$.
For given positive constants $u,T$, define the set of $r$th conjunctions
$ C_{r,T,u}:= \{t\in [0,T]: X_{r:n}(t) > u\}$
with $X_{r:n}(t)$ the $r$th largest order statistics  of $X_i(t),  t\ge 0, 1\le i\le n$.
In numerous applications such as brain mapping and digital communication systems, of interest is the approximation of the probability
that the set of conjunctions $C_{r,T,u}$ is not empty.
Imposing the Albin's conditions on $X$,
in this paper we obtain an exact asymptotic expansion of this probability as $u$ tends to infinity.
Furthermore, we establish the tail asymptotics of the supremum of  \cling{the order statistics processes of skew-Gaussian processes}
and
a Gumbel limit theorem for the {minimum} order statistics of stationary Gaussian processes.
As a by-product we derive a version of Li and Shao's normal comparison
lemma for the minimum and the maximum of Gaussian random vectors.
\end{abstract}

\section{Introduction and Main Result} \label{sec1}
Let $\{X(t),t\ge 0\}$ be a stationary process with almost surely (a.s.) continuous sample paths and denote by $X_1 \ldot X_n, n\inn $ \cling{mutually} independent copies of $X$.
\hH{Of} interest  in this contribution is the $r$th order statistics
process $X_{r:n}$
\tb{of} $X_1 \ldot X_n$, i.e., for any $t\ge 0$
\BQN \label{def: Orderp}
X_{n:n}(t)
 \le \cdots \le X_{1:n}(t).
\EQN
Throughout the paper, $X_{r:n}$ \eH{is referred to as} the $r$th order statistics process generated by the process $X$.
\eH{Order statistics play a central role in many statistical applications.
Naturally, the order statistics processes are of particular interest in statistical applications which involve the time-dynamics.}
\cL{If $X_i(t)$ is the value of a certain object (say image) $i$ \eH{measured at time point $t$},  and $u$ is a fixed threshold, then the set of points that the $r$th conjunction occurs \cling{before some time point $T$} is defined by}
$$ C_{r,T,u}:= \{t\in [0,T]: X_{r:n}(t) > u \}. $$
\eH{In applications it is of interest to calculate the probability that $C_{r,T,u}$ is not empty, which is given by}
\BQN\label{TailOS}
 p_{r,T}(u):=\pk{C_{r,T,u} \not=\phi}= \pk{\sup_{t\in [0,T]}X_{r:n}(t) > u}.
 \EQN
\cL{Clearly, in an engineering context where $X_i$'s model some random signals, $p_{r,T}(u)$ relates to the probability that
at least $r$ signals overshoot the threshold $u$ at some point during the time interval $[0, T]$. Most prominent statistical applications, concerned with the analysis of the surface roughness during all machinery processes and functional magnetic resonance imaging (\cling{FMRI}) data, relate to the calculations of $p_{r,T}(u)$.}
\eH{A methodology for the analysis of FMRI is established in the seminal contribution \cite{MR1747100}.
Therein the authors derive approximations of $p_{n,T}(u)$ by calculating the
expectation of Euler \tbb{characteristic} of $C_{n,T,u}$
for \cling{a fixed high threshold $u$}. }

For \hH{certain} smooth  Gaussian random fields approximations of $p_{n,T}(u)$ have been discussed in \cite{MR2775212,ChengXiao13,MR1747100},
whereas results for \hH{some} non-Gaussian random fields are derived in \cite{MR2654766}.
\tb{Exact asymptotic expansion of $p_{r,T}(u)$
for the class of stationary Gaussian processes $X$ was recently
derived in \cite{DebickiHJminima}.}
Obviously, the Gaussian random field cannot be used to model phenomena
and data sets that exhibit certain non-Gaussian characteristics such as skewness. \hH{It} arises in many \tb{applied-oriented} fields
\hH{including} engineering, medical \cling{science}, agriculture and environmental \cling{studies}; see, e.g., \cite{AlRawwashS07,AlodatR09,ZareifardJ2013}. In recent years, new technologies
such as FMRI and positron emission tomography have been used to collect data concerning the living human brain as well as astrophysics.
As mentioned in the literature, these images can be efficiently modeled by stationary random fields.

\cL{Since the exact calculation of $p_{r,T}(u)$ is not possible in general},
\eH{in this contribution} we \hH{derive} approximations of $p_{r,T}(u)$ for $u$ large.

For the formulation of our main result we need to introduce Albin's conditions \hH{imposed on $X$ as suggested in}  \cite{AlbinPHD,Albin1990,Albin2003}. In what follows, let
$\mathrm D$ be a non-empty subset of $\R$.
\\
{\bf Condition A($\mathrm D$)}: (Gumbel MDA \cling{and} conditional limit)
Suppose that $X(0)$ has a continuous df with infinite right endpoint, and it is in the Gumbel max-domain of attraction (MDA), i.e.,
for some positive scaling function $w(\cdot)$ we
have as $u\to \IF$
\BQN\label{Gumbel: c}
\pk{X(0)> u+ \frac x{w(u)}} = \pk{X(0)> u} e^{-x}(1+o(1)), \quad \forall x\inr.
\EQN
Let $q=q(u)$ \eH{satisfying} $\lim_{u\uparrow \IF} q(u) =0$ be a strictly positive non-increasing function. Assume that for any $y\in \mathrm D$ there exists a random process
$\{\xi_{y}(t), t\ge0\}$, such that
for any grid of points $0< t_1 < \cdots < t_d < \IF$
\eH{we have} the convergence in distribution (denoted by $\todis$)
\BQN \label{condL}
&&\nonumber\Bigl(w(u) (X(qt_1)-  u)  \ldot
w(u) (X(qt_d )-  u)\Bigr) \Big\lvert \Bigl( w(u) (X(0 )-  u) >y \Bigr)\\
&&\qquad
\todis \Bigl(\xi_{y}(t_1) \ldot \xi_{y}(t_d)\Bigr), \quad u\to \IF.
\EQN
{\bf Condition B}: (Short-lasting-exceedance) For all positive constants  $a, T$
\BQN\label{Cond: B}
\lim_{N\to \IF}\limsup_{u\to \IF} \sum_{k=N}^{[T/(aq)]} \pk{X(aqk)> u \lvert X(0) > u}  = 0,
\EQN
where $q$ is given as in condition A($\mathrm D$) and $[x]$ denotes the integer part of $x$.\\
{\bf Condition C}: Suppose that there exist positive constants
$\lambda_0, \rho, b, C$ and $d>1$  such that 
\BQN\label{Cond: C}
\pk{X(q t)> u+\frac\lambda{ w(u)},  X(0) \le u \Big \lvert X(q t) > u} \le C t^d\lambda^{-b}
\EQN
holds for all $u$ large and all  $t$ \hH{positive} such that $ 0 < t^\rho < \lambda < \lambda_0$. \cling{Here $w$ and $q$ are given as in condition  A($\mathrm D$)}.

Here we have chosen a simpler condition C than that in  \cite{AlbinPHD,Albin1990}.
It has been shown in \cite{Albin2003} that condition C  above is sufficient for the validity of
condition C given  in  \cite{AlbinPHD,Albin1990};
see also Proposition 2 in \cite{Albin1992}.

Note that the Albin's conditions A($\mathrm D$),B,C given above are satisfied by many well-known stationary processes \cL{such as $\hH{\chi^2}$,
$\Gamma$ and $\sqrt F$ processes in \cite{AlbinPHD}. \hH{A concrete example is  the Slepian process}, see \hH{for} other examples}  \cite{Albin2003,Aue2009,TurkmanA,TurkmanB}. However, showing the validity of these conditions requires in general significant efforts.

Let  $\xi^{(i)}_{0}, i\le  n$ be independent copies of \hH{the random process} $\xi_{0}$ given in condition A$(\{0\})$. In order to derive the exact asymptotics of $p_{r,T}(u)$ we introduce the following constants 
\BQN \label{Hr}
\mathcal A_{r}:= \lim_{a \downarrow 0} \frac{1}{a} \pk{ \sup_{k\ge 1}{\min_{1\le i\le r}\xi_0^{(i)}}(ak) \le 0 },\ \ r\le n,
\EQN
which we refer to as the \tn{{\it generalized Albin constants}}. The finiteness and positiveness of it will be established in \netheo{T1} below.  \hH{For notational simplicity we set hereafter
$ c_{r,n}={n!}/({r!(n-r)!})$.} \tn{Next, we state} our principle result.

\begin{theorem} \label{T1}
Let $\{X_{r:n}(t), t\ge 0\}$ be the $r$th order statistics process generated by the stationary process $X$.
\tb{If  \cling{conditions} A($\{0\}$), B and  C hold for $X$,}
then for any $T>0$, as $u\to \IF$,
\BQN
\pk{ \sup_{t\in [0,T]}
X_{r:n}(t) > u}&=& T \mathcal A_{r}\hH{c_{r,n}}\frac{\bigr( \pk{X(0)> u} \bigl)^r}{q(u)} (1+o(1))\eH{,} \
\EQN
where $\mathcal A_{r}$  defined in \eqref{Hr} is finite and positive.
\end{theorem}
This paper is organized as follows: In \nesec{sec2} we discuss an application of \netheo{T1}  concerning the skew-Gaussian processes
and then derive  the Gumbel limit theorem for the {minimum} order statistics process generated by a stationary Gaussian process $X$. 
All the proofs are presented in \nesec{sec3}. \nesec{sec5} gives an Appendix which establishes a version of \cling{Li and Shao\rq{}s normal comparison lemma} 
for {the minimum and maximum} order statistics of Gaussian random vectors.

\section{Skew-Gaussian Processes and Gumbel Limit Theorem}\label{sec2}
Throughout this section \tn{ we assume that  $\{X(t), t\ge0\}$ is} a centered stationary Gaussian process with a.s. continuous sample paths and covariance function $\rho(\cdot)$ \cL{such \hH{that}, for some $\alpha\in(0, 2]$}
\BQN\label{corrr}
\rho(t) < 1,\ \ \forall t>0 \ \ \text{and} \  \ \rho(t) = 1- \abs t^\alpha + o(\abs t^\alpha),
\ \ t\to0.
\EQN
It is known (\cling{see,} e.g., \cite{Albin1990,Pit96,Faletal2010}) that the process $X$ satisfies the assumptions of \netheo{T1} with the process $\xi_{0}$ in condition A$(\{0\})$ given by
\BQN \label{eq:eta}
\xi_{0}(t) = \sqrt 2Z(t) -t^\alpha + E, \quad t\ge0,
\EQN
where
\tb{$E$ is a unit exponential random variable (rv)
and
$\{Z(t), t\ge0\}$ is a (independent of $E$) standard fractional Brownian motion (fBm)
with Hurst index $\alpha/2 \in (0,1]$, i.e., $Z$ is a centered Gaussian
process with a.s.\ continuous sample paths and covariance function
}
$$
 \Cov(Z(s), Z(t))= \frac12 \Bigl( s^\alpha+ t^\alpha- \abs{s-t}^\alpha\Bigr), \quad s, t \ge 0.
$$
We note in passing that the findings of \netheo{T1} for such $X$ coincide with those of Theorem 2.2 in
\cite{DebickiHJminima}. Our setup here is however more general than that of the aforementioned paper,
\tn{ as we demonstrate \cling{below}.} 
Let $\{X_i(t),t\ge 0\},i\le m+1, m \inn$ be independent copies \hH{of $X$}.
For any   $\delta\in \hH{(0, 1]}$ define the  \tb{{\it skew-Gaussian}} process $\hH{\zeta}$ as
\BQN \label{def: SkewG}
\hH{\zeta}(t)
= \delta \abs{\X(t)}+\sqrt{1-\delta^2} X_{m+1}(t), \ \abs{\X(t)}=\big(\sum_{i=1}^m X_i^2(t)\big)^{\frac12}, \quad \hH{t\ge 0}\cling{.}\qquad
\EQN
\begin{theorem} \label{ThmA}
 Let $\{\hH{\zeta_{r:n}(t)}, t\ge0\}$ be the $r$th order statistics process generated by the skew-Gaussian process $\zeta$. 
\hH{If the stationary Gaussian process $X$ has covariance function $\rho(\cdot)$ which satisfies \eqref{corrr}},
then for any $T>0$, as $u\to\IF$,
\BQN\label{wideA}
&&\pk{ \sup_{t\in [0,T]} \cling{\zeta_{r:n}}(t) > u} \notag\\
&&\quad= T \widetilde{\mathcal{A}}_{r,\alpha} c_{r,n}\delta^{rm-r}\frac{2^{r-rm/2}}{(\Gamma(m/2))^r} u^{\frac2\alpha +rm-2r}
e^{-\frac{ru^2}2}(1+o(1)),
\EQN
where $\widetilde{\mathcal{A}}_{r,\alpha}$ is determined by \eqref{Hr}
with $\xi_0^{(i)}, i\le n$ being $n$ independent copies of $\xi_{0}$ given in \eqref{eq:eta}, \cling{and $\Gamma(\cdot)$ stands for the Euler Gamma function}.
\end{theorem}
\textbf{Remarks}.
\eH{a)} The \tn{special case of \netheo{ThmA}, for $\delta=1$,} coincides with that obtained in Corollary 7.3 in \cite{Pit96}; see also \cite{Faletal2010,HJ14}.\\
\eH{b)} The \eH{Pickands constant $\mathcal{H}_\alpha$ coincides with $\widetilde{\mathcal{A}}_{1,\alpha}$ if $n=1$. It is well-known that
$\mathcal{H}_{1}=1$ and $\mathcal{H}_2= 1 /\sqrt{\pi}$. For other values of $\alpha$ the recent contribution \cite{DiekerY} (see also the excellent monograph \cite{Yakir}) suggests an efficient algorithm to  simulate  $\mathcal{H}_\alpha$. For $n>1$ both calculation and simulation of $\widetilde{\mathcal{A}}_{r,\alpha}$ are open problems.}\\

In extreme value analysis (see, e.g., \cite{Berman92,leadbetter1983extremes,Faletal2010}) it is also of interest
to find some normalizing constants $a_T >0, b_T\in \R$
so that the linear normalization of the
supremum $a_T \Big(\sup_{ t\in [0, T]} X_{r:n}(t)- b_T\Big)$ converges in distribution as $T\to\IF$, where $X_{r:n}$ is the $r$th order statistics process generated by the stationary Gaussian process $X$.
The following theorem gives a Gumbel limit result for $X_{n:n}$ {generated by}
\hH{a}  weakly dependent stationary Gaussian process.

\begin{theorem}\label{ThmB}
Let $\{X_{n:n}(t),t\ge0\}$ be the minimum order statistics process generated by the stationary Gaussian process $X$ with covariance function $\rho(\cdot)$ satisfying \eqref{corrr}.
If  $\rho(t) \ln t= o(1) $ \eH{holds} as $t\to\IF$, then
\BQN\label{GuG}
 \limit{T}\sup_{x\in \R} \Abs{ \pk{ a_T \Big(\sup_{ t\in [0, T]} X_{n:n}(t)- b_T\Big) \le x} -
\expon{ - e^{-x}}} =0,
\EQN
where (set below  $D:= (n/2)^{n/2 -1/ \alpha}\widetilde{\mathcal{A}}_{n,\alpha}(2\pi)^{-n/2}$) 
\BQN \label{Normalized constant}
 a_T = \sqrt{2n\ln T}, \  b_T =\sqrt{\frac{2\ln T}n} + \frac{
\left(\frac 1\alpha -\frac n2\right)\ln\ln T + \ln D}{\sqrt{2n\ln T}}.
 \EQN
\end{theorem}
\textbf{Remarks}. a) It follows from the proof of \netheo{ThmB} that
a similar result still holds for the maximum order statistics process $X_{1:n}$ under the same condition
(\tb{since} \eqref{eq:Min} holds for the maximum).\\
b) In several applications it is of interest to consider a random time interval \cling{$\mathcal{T}$} 
instead of $T$; see\cling{,} e.g.,
\cite{FERET,TANWU}. As in \cite{TANWU} our result in \eqref{GuG} can be extended for random intervals; we omit that result since it can be shown with
similar arguments as in the aforementioned paper. \\
c) The deep contribution \cite{Oleg05} shows that besides Gumbel limit theorems, of interest for \cL{applications} is
the growth of $\E{(\sup_{ t\in [0, T]} \abs{X_{n:n}(t)}) ^p}$ for given $p>0$. 
In view of \netheo{ThmA} (with $m=\delta=1$) and applying Lemma 4.5 in \cite{TANSPA} we obtain 
\BQN
 \limit{T} \E{ \Bigl( \frac{ \sup_{ t\in [0, T]} \abs{X_{n:n}(t)}}{\sqrt{2/n}\ln T} \Bigr)^p}= 1.
\EQN

\section{Proofs}\label{sec3}
In this section, we present \tb{proofs} of Theorems~\ref{T1}, \ref{ThmA} and \ref{ThmB}.
We shall rely on the methodology developed in the seminal paper \cite{Albin1990}. As mentioned therein, checking the Albin's conditions
for stationary processes is usually a hard task.
\tb{In \cling{Section}~\ref{proof.th1} we consider} $X$ to be a stationary process with a.s.
continuous sample paths. \tb{In Sections~\ref{proof.th2}, ~\ref{proof.th3}}
we concentrate on the special case where $X$ is a centered stationary Gaussian
process with a.s. continuous sample paths and covariance function
$\rho(\cdot)$ satisfying \eqref{corrr}.

\subsection{Proof of Theorem \ref{T1}}\label{proof.th1}

We \tb{begin with some preliminary lemmas that will be used in the proof of \netheo{T1}.}
\eH{Unless otherwise specified,  $\{X_{r:n}(t), t\ge0\}$ denotes the $r$th order statistics process generated by the stationary process $X$.}

\tb{The next} lemma plays a key role throughout the proofs.
\tb{Since its proof is straightforward, we omit it. \hH{Recall that we defined} $c_{r,n}={n!}/({r!(n-r)!})$.}
\begin{lemma}
\eH{\hH{If $X(0)$ has a continuous distribution \cling{function}, then} for any} $t\ge0$
\BQN \label{eq:Xrnt}
\pk{X_{r:n}(t) > u} = c_{r,n} \bigr(\pk{X(t)>u}\bigl)^r (1+o(1)), \quad u\to \IF.
\EQN
\end{lemma}
\begin{lemma} \label{AlbA}
If condition A($\mathrm D$) \tb{holds for $X$}, then
$X_{r:n}(0)$ has df  in the Gumbel MDA with scaling function $w_r(u) = rw(u)$.
Further, for any grid of points $ 0<t_1< \cdots<t_d<\IF$  and all $y\in \mathrm D$
we have 
\BQN \label{eq:ConA}
&&\Bigl( w_r(u)(X_{r:n}(qt_1)-u) \ldot
w_r(u)(X_{r:n}(qt_d)-u) \Bigr) \Bigl \lvert \Bigl(w_r(u)(X_{r:n}(0)-u)> {r}y \Bigr)
 \nonumber\\
&& \qquad \todis
\Bigl( \min_{1 \le i \le r}  r \xi_{y}^{(i)}(t_1)   \ldot \min_{1 \le i \le r}  r \xi^{(i)}_{y}(t_d) \Bigr), \quad u\to\IF,\qquad\qquad
\EQN
where \tb{$\xi_{y}^{(i)}, i\le n$ are mutually independent copies of
$\xi_{y}$ as in condition A($\mathrm D$).}
\end{lemma}
\textbf{Proof}.
First,   by \eqref{Gumbel: c} and \eqref{eq:Xrnt}
\BQNY
\pk{ X_{r:n}(0) > u + \frac{x}{r w(u)} } =\pk{X_{r:n}(0) >u} e^{-x} (1+o(1)) , \quad x\in \R
\EQNY
meaning that $ X_{r:n}(0)$ has df in the Gumbel MDA with $w_r(u)=r w(u)$.\\
Further, it follows from \eqref{condL} that the  convergence in distribution
\BQN\label{eii}
\biggl( X_{i}^*(qt_1
) \ldot X_{i}^*(qt_d) \Bigr) \biggl \lvert (X_{i}^*(0)>{r}y) \todis
\bigl( r\xi^{(i)}_{y}(t_1)\ldot r \xi^{(i)}_{y}(t_d)\bigr)
\EQN
holds \hH{as $u\to\IF$} for all $i\le n, y\in \mathrm D$, where $\cL{X_{i}^*(t)}= w_r(u)
(X_{i}( t)- u)$. Let further $\cL{Y_{r}^*}(t)= w_r(u)\bigr(X_{r:n}(t)-u\bigl)$ and fix a  grid of points
$0<t_1< \cdots<t_d<\IF$. Next, we show that \eqref{eq:ConA} holds when $r=n$.
\tb{Indeed}, for any given constants $y_1 \ldot y_d \in\R$ \hH{by   \eqref{eii}}
\BQN \label{decomp: Min}
&&\pk{ Y_n^*(qt_1)> y_1 \ldot
Y_n^*(qt_d)> y_d \Bigl \lvert Y_n^*(0)> ny }
 \notag\\
&& = \frac{\pk{ \min_{1 \le i \le n} X_{i}^*(qt_j)> y_j, 1\le j\le d , \min_{1 \le i \le n}
X_{i}^*(0)>ny} } {\pk{\min_{1 \le i \le n}X_{i}^*(0)> ny}}
\notag\\
&& \to \pk{\min_{1 \le i \le n}  n \xi^{(i)}_{y}(t_1) > y_1  \ldot
\min_{1 \le i \le n}  n \xi_{y}^{(i)}(t_d) > y_d}
\EQN
as $u\to\IF$. Similarly, the claim of \eqref{eq:ConA} holds for all $r<n$ if we show that,
for any given constants $y_1 \ldot y_d \in\R$
\BQN \label{decomp: A}
&&
\pk{ Y_r^*(qt_1)> y_1 \ldot
Y_r^*(qt_d)> y_d \Bigl \lvert {Y_r^*(0)> ry}}
 \notag\\
&& = \frac{ \pk{ \min_{1 \le i \le r} X_{i}^*(qt_j)> y_j, \cling{1\le j\le d}, \min_{1 \le i \le r}
X_{i}^*(0)>ry} } {\pk{\min_{1 \le i \le r}X_{i}^*(0)>
ry}}\notag\\
&&\quad\times(1+ \Upsilon_r(u)),\ \mbox{with}\ \lim_{u\to\IF} \Upsilon_r(u)=0.
\EQN
\hH{Next}, we  only present the proof for the case that $r= n-1$ and
 \tb{$d=1$; the other cases follow by} similar arguments. \hH{By} \eqref{eq:Xrnt}
 $$\pk{Y_{n-1}^*(0) > (n-1)y} = n \pk{\min_{1\le i\le n-1} X_{i}^*(0) >(n-1)y} (1+o(1)) $$
 as $u\to\IF$. \hH{Further}
\BQNY
 \lefteqn{ \pk{Y_{n-1}^*(qt_1) > y_1, Y_{n-1}^*(0) >(n-1)y}}
 \\
 &= & \pk{Y_{n-1}^*(qt_1) > y_1 \ge Y_n^*(qt_1), Y_n^*(0) > (n-1)y}
  \\
 &&  + \pk{Y_n^*(qt_1) > y_1, Y_{n-1}^*(0) > (n-1)y \ge Y_n^*(0)}\\
 &&  + \pk{Y_{n-1}^*(qt_1) > y_1 \ge Y_n^*(qt_1), Y_{n-1}^*(0) > (n-1)y \ge
 Y_n^*(0)} \\
 && +  \pk{Y_n^*(qt_1) > y_1, Y_n^*(0) > (n-1)y}\\
 & =:& I_{1u} + I_{2u} + I_{3u} + {I_{4u}}.
  \EQNY
  Since as $u\to\IF$
\BQNY
&&\pk{X_n^*(qt_1) \le y_1, X_n^*(0) > (n-1)y}\cling{\le  \pk{X_n(0) > u+\frac y{w(u)}} =o(1)}, 
\EQNY
 we have
  \BQNY
 I_{1u} &=& 
   n \pk{ \min_{1 \le i \le n-1} X_{i}^*(qt_1)>
 y_1, \min_{1 \le i \le n-1} X_{i}^*(0)>(n-1)y } o(1),
  \EQNY
 and similarly $I_{2u} =  I_{1u}(1+o(1)).$
Using further the fact that
$$\pk{X_n^*(qt_1) \le y_1, X_n^*(0) \le
(n-1)y} = 1+o(1), \ u\to\IF,$$
 we have
 \BQNY
I_{3u}&=& \sum_{i, i'\le n}\mathbb P\left(\min_{1\le j \le n, j\neq i} X_{j}^*(qt_1) >y_1, X_{i}^*(qt_1) \le y_1,\right.\quad\\
&&\quad\left.\min_{1\le j' \le n, j'\neq i'} X_{j'}^*(0) >(n-1)y, X_{i}^*(0) \le (n-1)y\right)
 \quad \\
 & =& n \pk{\min_{1\le j \le n-1} X_{j}^*(qt_1) >y_1, \min_{1\le j' \le n-1} X_{j'}^*(0) >(n-1)y}\\
 &&\times \pk{X_n^*(qt_1) \le y_1, X_n^*(0) \le (n-1)y} \\
 &&  + \cling{c_{2,n}}\pk{\min_{1\le j \le n-2} X_{j}^*(qt_1) >y_1, \min_{1\le j' \le n-2} X_{j'}^*(0)
 >(n-1)y} \\
 &&\times  \pk{{X_{n-1}^*}(qt_1) \le y_1, X_{n-1}^*(0) > (n-1)y}\pk{X_n^*(qt_1) > y_1, X_n^*(0) \le
 (n-1)y} \\
 & =& n \pk{\min_{1\le j \le n-1} X_{j}^*(qt_1) >y_1, \min_{1\le j' \le n-1} X_{j'}^*(0)
 >(n-1)y} (1+ o(1)).
\EQNY
Since in view of \eqref{decomp: Min}, for $k=0, 1, 2$, \hH{as $u\to\IF$},
\BQNY
\pk{\min_{1\le j \le n-k}
X_{j}^*(qt_1) >y_1, \min_{1\le j' \le n-k} X_{j'}^*(0)
 >(n-1)y} = \bigr(\pk{X(0) >u}\bigl)^{n-k} O(1),
\EQNY
we conclude that \cling{$I_{4u}=I_{3u}o(1)$, and further that} \eqref{decomp: A} holds for  $r= n-1$ and
 $d=1$.
This completes the proof.  \qed 
\begin{lemma} \label{AlbB}
If condition B is satisfied by \tb{$X$}, then for any $a,T$ positive
 \BQN\label{Cond: B}
\lim_{N\to \IF}\limsup_{u\to \IF} \sum_{k=N}^{[T/(aq)]}
\pk{X_{r:n}(aqk)> u \lvert X_{r:n}(0) > u} = 0.
\EQN
\end{lemma}
\textbf{Proof}. First, since
for all integers $k\ge 1$ and any $u$ positive
\BQNY
\lefteqn{
\pk{X_{n:n}(aqk)> u \lvert X_{n:n}(0) > u}
 = \frac{\pk{X_{n:n}(aqk) >u, X_{n:n}(0)>u }}{
\pk{X_{n:n}(0)>u}} } \\
&& = \Bigr(\pk{X(aqk) > u \lvert X(0)> u}\Bigl)^n \le \pk{X(aqk) > u
\lvert X(0)> u}
\EQNY
holds, \hH{condition B  implies the claim for $r=n$.} 
\hH{If $r<n$, then with similar arguments as in \eqref{decomp: A} we have} \cling{for large $u$}
\BQNY
\pk{X_{r:n}(aqk) >u\Big\lvert X_{r:n}(0) >u}
&=& \left(\pk{X(aqk)> u\Big\lvert X(0) >u}\right)^r (1+\Upsilon_r(u))
\\
&\le&
 K \Bigr(
\pk{X(aqk)> u \Big \lvert X(0)> u} \Bigl)^r
\EQNY
holds for some $K>0$, \hH{hence again condition B establishes the proof.} \qed
\begin{lemma} \label{AlbC}
If condition C is satisfied by \tb{$X$} with the parameters as therein,
then there exists some  positive constant $C^*$ such that for all $u$ large 
\BQN\label{lemCC}
\pk{X_{r:n}(q t)> u+
\frac\lambda{ w(u)}, X_{r:n}(0) \le u \Big \lvert X_{r:n}(q t) > u}
\le  C^* t^d\lambda^{-b}
\EQN
holds for any $t$ \hH{positive} satisfying $ 0 < t^\rho < \lambda < \lambda_0$.
\end{lemma}
\textbf{Proof}. 
\hH{By} condition C,  for \tn{sufficiently} large $u$ and $C^*= nC, \hH{r=n}$
\BQN \label{Cond: Cmin}
\lefteqn{\pk{X_{n:n}(q t) >u+\frac\lambda{ w(u)},  X_{n:n}(0) \le u \Big \lvert X_{n:n}(q t) > u} } \notag \\
&& \le \sum_{i=1}^n \frac{\pk{X_{n:n}(q t) > u+\frac\lambda{ w(u)},
X_i(qt) \le u}}
{\pk{ X_{n:n}(qt)>u}} \notag \\
&& \le  \sum_{i=1}^n \frac{\pk{X_i(q t) > u+\frac\lambda{ w(u)},
X_i(0) \le u, \min_{1 \le j \le n, j\not=i} X_j(q t) > u} }{(\pk{ X(qt)> u})^n} \notag \\
&& =   \sum_{i=1}^n \frac{\pk{X_i(q t) > u+\frac\lambda{ w(u)},  X_i(0) \le u} }{\pk{ X(\cling{qt})> u}} \notag \\
&& =   \sum_{i=1}^n \pk{X_i(q t)  > u+\frac\lambda{ w(u)},  X_i(0) \le u \Big \lvert   X_i(q t)  > u} \notag  \\
&& \le   C^* t^d\lambda^{-b}
\EQN
\cling{holds} for all  $t$ \cling{positive} satisfying $0< t^\rho < \lambda
< \lambda_0$.
\hH{If $r<n$, then with similar arguments}
as in \eqref{decomp: A}
\BQNY \lefteqn{
\pk{ X_{r:n}(qt)
> u+ \frac\lambda{w(u)}, X_{r:n}(0) \le u}
} \\
&& =  c_{r,n}\pk{ \min_{1\le i\le
r}X_i(qt)>u+\frac\lambda{w(u)}, \min_{1\le i\le r}X_i(0) \le u}(1+o(1))
\EQNY
\hH{holds} as $u\to\IF$. Consequently, it follows from \eqref{Cond: Cmin} that  there exists some positive constant $C^*$ such that
\BQNY
\pk{X_{r:n}(q t)>u+ \frac\lambda{ w(u)}, X_{r:n}(0) \le \cling{u}\Big \lvert X_{r:n}(q t) > u}
\le  C^* t^d\lambda^{-b}
\EQNY
holds for all $t>0$ and $0< t^\rho < \lambda < \lambda_0$ \hH{establishing thus the proof}\cling{.}
\qed

\noindent\textbf{Proof of Theorem}~\ref{T1}.
The proof is based on an application of Theorem 1 in \cite{Albin1990}, see also Lemma A in \cite{Albin2003}. It follows from   Lemmas  \ref{AlbA}, \ref{AlbB} and \ref{AlbC}  that the conditions of \hH{Theorem 1 in \cite{Albin1990} are satisfied}, hence
for any $T>0$
\BQNY
\pk{\sup_{t\in[0, T]} X_{r:n}(t) >u} = T \mathcal A_{r} \frac{\pk{X_{r:n}(0) >u }}{q(u)} (1+o(1)), \quad u\to \IF,
\EQNY
where \hH{the \cling{limit}} in the right-hand side of \eqref{Hr} \hH{exists \cling{with $\mathcal A_{r}\in(0, \IF)$}.}
\hH{Hence the proof follows from \eqref{eq:Xrnt}.} \qed

\subsection{Proof of Theorem \ref{ThmA}}\label{proof.th2}
In the following, we \tb{focus on} the special case where the process $X$ is a
centered stationary Gaussian process with a.s. continuous sample paths and covariance function $\rho(\cdot)$ satisfying \eqref{corrr}.
\tb{Before proceeding to the proof of \netheo{ThmA}, we present \hH{four} lemmas.}
\begin{lemma} \label{LTA1}
\hH{If}  $\{\zeta (t), t\ge0\}$ is given as in \eqref{def: SkewG}, then 
\BQNY
\pk{\zeta(0) > u} =
\delta^{m-1}  \frac{2^{1-m/2}}{\Gamma(m/2)} u^{m-2} \expon{-\frac{u^2}2} (1+o(1))
\EQNY
and \hH{further the convergence in probability}
\BQN \label{eq: CondIF}
\tb{\abs{\X(0)} \Big\lvert \bigl(\zeta (0)>u \bigr) \toprob \IF}
\EQN
holds as $u\to \IF. $
\end{lemma}
\textbf{Proof}.
\hH{Since} $\abs{\X(0)}^2/2$ has  $Gamma(m/2, 1)$ distribution, \hH{we have }
\BQNY
\pk{\abs{\X(0)} >u} = \frac{2^{1-m/2}}{\Gamma(m/2)} u^{m-2} \expon{-\frac{u^2}2} (1+o(1)), \quad u\to\IF.
\EQNY
\hH{Hence}, Theorem 2.2 in \cite{FarkasH13} implies for any $\delta\in(0, 1]$
\BQNY
\pk{\zeta(0) > u} = \delta^{m-1} \hH{\pk{\abs{\X(0)} >u}} (1+o(1)), \quad u\to \IF.
\EQNY
Clearly, \eqref{eq: CondIF} holds for $\delta=1$. Next, taking a constant $c$ such that $1<c < 1/\sqrt{1-\delta^2}$ for $\delta\in(0,1)$, it
follows from the proof of Lemma 2.3 in \cite{FarkasH13} that
\BQNY
\hH{\limit{u}}\pk{(1-c \sqrt{1-\delta^2})u < \delta\abs{\X(0)}< c\delta u \Big \lvert \zeta(0) > u} =1,
\EQNY
implying thus \eqref{eq: CondIF}. \tbb{Hence the proof is complete}.
\qed

\begin{lemma} \label{LTA2}
Let $\{\zeta(t), t\ge0\}$ be given as in \eqref{def: SkewG}. \hH{If}
the covariance function $\rho(\cdot)$ of the generic stationary Gaussian process $X$
satisfies \eqref{corrr}, then for any grid of points
$0< t_1 < \cdots < t_d<\IF$  the joint convergence in distribution
\BQNY
\Bigl( u( \zeta (qt_1)- u) \ldot
 u ( \zeta (qt_d)- u) \Bigr) \Bigl \lvert ( \zeta (0)>u)
\todis  \Bigl( \xi_0(t_1) \ldot \xi_0(t_d)\Bigr)
\EQNY
holds as $u\to \IF$, where the process $\xi_0 $ is given by \eqref{eq:eta} and
\hH{ $q\cling{=q(u)}:=u^{-2/\alpha}$.}
\end{lemma}
\textbf{Proof}.
By \nelem{LTA1}, for any $s\inr$
$$
 \limit{u} \frac{\pk{ \zeta (0)> u + \frac s{u}}}{\pk{   \zeta (0)> u}} = e^{-s}.
$$
Consequently, we have the convergence in distribution
$$
 u(  \zeta (0)- u) \Bigl\lvert  (  \zeta (0)- u>0)\ {\todis}\   E, \quad u\to \IF,
 $$
with $E$ \hH{a} unit exponential rv. In view of Theorem 5.1 in \cite{Berman82}, it suffices to show that as $u\to\IF$
 \BQNY
 && \Bigl( u( \zeta (q t_1)- u ) \ldot u( \zeta (qt_d)- u) \Bigr) \Bigl
\lvert (  \zeta (0)= {u_x})\nonumber\\
 && \ \ \ \  \todis
\Bigl( \sqrt 2Z(t_1) - t_1^\alpha +x \ldot \sqrt 2Z(t_d) - t_d^\alpha +x\Bigr) , \quad u_x:= u+ x/u
\EQNY
 holds for all $d\ge1$ and almost all \hH{$x>0$ with 
 $Z$ a standard fBm \cling{with Hurst index $\alpha/2$}.}
Define
$$
X_{i}^*= \rho(q t_j)X_i(0), \quad
\cling{\Delta_i}(t_j)= X_i(q t_j) - X_{i}^*, \ \ \ i\le m+1,\ j\le d.
$$
For any $u>0$ and $j\le d$, we have
\begin{align*}
& {\big(u[ \zeta (q t_j)- u]  -x\big) \Bigl \lvert (  \zeta (0)=u_x) }
\\
 \nonumber 
& = \cling{\Bigg(} \delta u\left(\sqrt{\sum_{i=1}^m  \big(X_i^2(0)  +  2 X_{i}^*  \Delta_i(t_j) - (1-\rho^2(qt_j))X^2_{i}(0) +
\Delta^2_i(t_j)\big)
}  - \hH{\abs{\X(0)}} \right)
\\
\nonumber & \quad   +  \sqrt{1-\delta^2}u
\Bigl(\Delta_{m+1}(t_j) - (1-\rho(qt_j))X_{m+1}(0) \Bigr)\Bigg)\Biggl \lvert
(  \zeta (0)=u_x)
\\
 \nonumber & =:  \delta A_u  + \sqrt{1-\delta^2} B_u .
 \end{align*}
Let $ Z_i,   i\le m+1 $ \hH{be \cling{mutually} independent copies of $Z$.}
In view of \eqref{corrr}  for $s,t >0$ and $i\le  m+1$
\BQN\label{DeltaA}
\hH{\limit{u}}u^2 \Cov(\Delta_i(s),\Delta_i(t)) =   s^{\alpha} +
t^{\alpha} - \abs{s-t}^\alpha = 2\Cov(Z_i(s), Z_i(t)), \quad  \quad
\EQN
which implies  the \hH{following} convergence of finite-dimensional \tb{distributions}
\BQNY
\{u\Delta_i(t), t\ge 0\}\todis\{\sqrt{2}Z_i(t),t\ge 0\},\ \ u \to\IF,\ \ i\le  m+1.
\EQNY
By the independence of $\Delta_i$'s 
 and \cling{$X_{i}$'s}, the $Z_i$'s can be chosen \hH{to be} independent of $ \zeta (0)$.
Further, since $(X_1(0) \ldot X_{m+1}(0))$ is a centered Gaussian random
 vector with $N(0,1)$ independent components, we have the stochastic representation (\cling{see} \cite{Cambanis81})
\BQN \label{Exp1}
(X_1(0) \ldot X_m(0), X_{m+1}(0)) \equivdis R(\vk{O}B, I \sqrt{1-B^2}),
\EQN
where $\vk{O}= (O_{1}, \ldots, O_{m})$ is a random vector
uniformly distributed on the unit sphere of $\R^m$. Here the rv  $I$ satisfies
$\pk{I=\pm1} = 1/2$, the random radius $R>0$ a.s. is such that $R^2$ has chi-square
distribution with $m+1$ degrees of freedom, and  the rv $B$ is \tb{supported} in $(0,1)$ a.s.\ such that
$B^2$ has \hH{beta distribution} with parameters $m/2, 1/2$.
Moreover, $\vk{O}, I, R, B, Z_i, i\le m+1$ are mutually independent.
Consequently, using 
the fact that \cling{$\sqrt{x_0+x}= \sqrt {x_0}+ ({2\sqrt {x_0}})^{-1} x(1+o(1))$ as $x\to0$}, together with
\eqref{DeltaA} and \eqref{Exp1}, we obtain as $u\to\IF$
\begin{align*}
A_u & \cling{= \frac{\sum_{i=1}^m \rho(qt_j)X_i(0) [u\Delta_i(t_j)]- \sum_{i=1}^m\frac{u(1-\rho^2(qt_j))}{2}  X^2_i(0)+ \sum_{i=1}^m \frac{u\Delta^2_i(t_j)}2} {\abs{\X(0)}}} \\
&\quad \cling{\times(1+o_p(1)) \Bigl \lvert (  \zeta (0)= {u_x})}\\
& = \frac{\sum_{i=1}^m \sqrt 2X_i(0)Z_i(t_j) - \frac{ \sum_{i=1}^m X^2_i(0)}{u} t_j^{\alpha}+ \frac{\sum_{i=1}^m Z^2_i(t_j)}{u} }{\abs{\X(0)}} (1+o_p(1)) \Bigl \lvert (  \zeta (0)= {u_x})\\
&  \equivdis  \left(\sqrt 2\sum_{i=1}^mO_iZ_i(t_j) - \frac{ RB}{u} t_j^{\alpha}+ \frac{\sum_{i=1}^m Z^2_i(t_j)}{u RB} \right) \\
&\quad\times  {(1+o_p(1))} \Bigl \lvert (R(\delta B +\sqrt{1-\delta^2} \sqrt{1-B^2} I)=u_x)\\
B_u & \equivdis \Bigl(\sqrt 2 Z_{m+1}(t_j) -
\frac{R\sqrt{1-B^2}I}{u} t_j^{\alpha} \Bigr) (1+o_p(1)) \Bigl \lvert
(  \zeta (0)=u_x).
\end{align*}
Since the following stochastic representation 
$$\sum_{i=1}^mO_iZ_i(t_j) \equivdis Z_1(t_j)\Big(\sum_{i=1}^m O_i^{\Ae{2}}\Big)^{1/2} {= Z_1(t_j)}
$$
holds, we have further  by \eqref{eq: CondIF}
\BQNY
 &&\delta A_u  + \sqrt{1-\delta^2} B_u \\
&& \equivdis \left( \sqrt 2\Bigl(\delta\sum_{i=1}^mO_iZ_i(t_j)+\sqrt{1-\delta^2} Z_{m+1}(t_j)\Bigr) - \frac{R(\delta B + \sqrt{1-\delta^2} \sqrt{1-B^2}I)}{u}t_j^\alpha \right. \\
&& \quad \left. + \frac{\sum_{i=1}^m Z^2_i(t_j)}{u RB} \right)(1+o_p(1))  \Bigl \lvert (R(\delta B + \sqrt{1-\delta^2} \sqrt{1-B^2}I) =u_x) \\
&& \todis \sqrt 2Z(t_j) - t_j^\alpha, \quad u\to\IF
\EQNY
establishing the convergence for any fixed $t_j>0$.
The joint convergence in distribution for $0<t_1< \cdots <t_d<\IF$ can be shown with
similar arguments and is therefore omitted here. \qed

\begin{lemma} \label{LTB}
\COM{Let $\{ \zeta (t), t\ge0\}$, $\delta\in(0, 1]$  be given as in \eqref{def: SkewG},
where the covariance function $\rho(\cdot)$ of the generic stationary Gaussian
process $X$ satisfies \eqref{corrr}.
Then, with $q=q(u)=u^{-2/\alpha}$,}
\hH{Under the assumptions and the notation of \nelem{LTA2},   for any $a,T$ positive}
\BQNY
\limsup_{u \to \IF} \sum_{j= N}^{[T/(aq)]} \pk{ \zeta (aqj) > u\Big \lvert   \zeta (0) > u } \to 0, \quad N\to \IF.
\EQNY
\end{lemma}
\textbf{Proof}.
It follows from \eqref{corrr} that for any $\ve>0$ small enough
$$
 \frac12 t^\alpha  \le 1- \rho(t) \le 2 t^\alpha, \quad \forall t\in (0, \epsilon].
$$
Denote by $\X(qt) - \rho(qt) \X(0)  =(X_1(qt) - \rho(t) X_1(0), \ldots, X_m(qt) - \rho(qt) X_m(0))$, and define
\BQN\label{def: zeta*}
 \cling{\zeta^*}(qt) \equiv\delta\abs{\X(qt) - \rho(qt) \X(0)} + \sqrt{1-\delta^2} ( X_{m+1}(qt) - \rho(qt) X_{m+1}(0)).\qquad
\EQN
Since $X(qt) - \rho(qt) X(0)$ is independent of $X(0)$, and $X(qt) - \rho(qt) X(0) \equivdis\sqrt{1-\rho^2(qt)} X(0)$, we have that $\zeta^*(qt) $ is independent of $ \zeta (0)$, and
$$
\zeta^*(qt) \equivdis \sqrt{1-\rho^2(qt)} \zeta (0).
$$
Moreover, by the triangle inequality
\BQNY
\zeta^*(qt) &\ge &\Bigl(\delta \abs{\X(qt)} + \sqrt{1-\delta^2} X_{m+1}(qt)\Bigr)- \rho(qt) \Bigl(\delta \abs{\X(0)} + \sqrt{1-\delta^2} X_{m+1}(0)\Bigr)
\\
& =&  \zeta (qt) - \rho(qt) \zeta (0) > u (1-\rho(qt))
\EQNY
provided that $ \zeta (qt) >  \zeta (0) >u.$ Therefore
\BQNY
 \pk{ \zeta (qt) >u\Big\lvert  \zeta (0)>u}
& \le&
 2 \frac{\pk{ \zeta (qt)>  \zeta (0)>u}}{\pk{  \zeta (0)>u}}
 \\
& \le& 2 \frac{\pk{\zeta^*(qt) > u(1-\rho(qt)),  \zeta (0)>u}} {\pk{  \zeta (0)>u}}
\\
& =& 2 \pk{ \zeta (0) > u\sqrt{\frac{1-\rho(qt)}{1+\rho(qt)}}}.
\EQNY
Furthermore, it follows from Chebyshev's inequality and
\nelem{LTA1} \tn{that} for any $p > m$
\begin{align}  \nonumber
\pk{ \zeta (qt) > u \Big \lvert   \zeta (0) > u}
&\le
\left\{
\begin{array}{ll}
  2^{1+p} \frac{\E{\abs{ \zeta (0)} ^p}}{t^{\alpha p/2}}, & qt\in (0, \epsilon], \\
  2 \pk{ \zeta (0) > u\sqrt{\frac\lambda2}}, & qt\in (\epsilon, T]
\end{array}
\right. \\
\label{eq:shortLE} &\le
\left\{
\begin{array}{ll}
  K_p t^{-\alpha p/2}, & qt\in (0, \epsilon], \\
  K_p u^{m-1-p}, & qt\in (\epsilon, T]
\end{array}
\right.
\end{align}
is satisfied for some positive constant $K_p$, where {$\lambda =
1-\sup_{\epsilon <s \le T} \rho(s) >0$}, and the second inequality is
due to the fact that
$$\pk{ \zeta (0) >u} \le
C u^{m-1}\frac1{\sqrt{2\pi}u}\expon{-\frac{u^2}2}\le C_p u^{m-1-p}, \quad u>0
$$
holds for some positive constants $C$ and $C_p$. \hH{Hence, if}  $p=2(2/\alpha+m-1)$, then  for $t\ge1$
\BQNY
\pk{ \zeta (qt) > u \Big
\lvert  \zeta (0) > u} &\le& K_p( 1+ T^{\alpha(p-m+1)/2}) \max( t^{-\alpha
p/2}, t^{-\alpha(p-m+1)/2}) \\
& \le & C_p t^{-2}, \quad qt \in(0, T].
\EQNY
Consequently,
 \BQNY
\limsup_{u\to\IF} \sum_{j = N}^{[T/(aq)]}\pk{ \zeta (aqj) > u \Big \lvert
 \zeta (0)>u}\le C_p\sum_{j=N}^\IF (aj)^{-2} \to 0,\ \ N\to\IF
 \EQNY
establishing the proof. \qed

\begin{lemma} \label{LTC}
\COM{Let $\{ \zeta (t), t\ge0\}$, $\delta\in(0, 1]$  be given as in \eqref{def: SkewG},
where the covariance function $\rho(\cdot)$ of the generic stationary Gaussian process $X$
satisfies \eqref{corrr}.
Then, with $q=q(u)=u^{-2/\alpha}$,}
\hH{Under the assumptions and the notation of \nelem{LTA2}} there exist positive constants $C, p,\lambda_0, u_0$ and $d>1$ such that
\BQNY
\pk{ \zeta (qt) > u +\frac{\lambda}{u},
 \zeta (0) \le u} \le Ct^d\lambda^{-p} \pk{ \zeta (0)>u}
\EQNY
 for any \hH{positive} $t$ satisfying $0< t^{\alpha/2} < \lambda < \lambda_0$ and \hH{all} $u> u_0$.
\end{lemma}
\textbf{Proof}. \hH{By \eqref{corrr} there exists} $\epsilon >0$ such that
\BQNY \rho(t) \ge \frac12
\quad\mbox{and}\quad 1-\rho(t) \le 2t^\alpha
\EQNY
for all $t\in (0, \epsilon]$. Further, for any $t$ \hH{positive} satisfying $0 < t^{\alpha/2} < \lambda <
\lambda_0 := \min(1/8, \epsilon^{\alpha/2})$ and $u > 1$
$$
\frac{1}{\rho(qt)} - 1 \le 4 \frac{t^\alpha}{u^2} \le
\frac{\lambda}{2u^2}.
$$
Next, for
$$
\X_{1/\rho}(qt)  = (X_1(qt) - \rho^{-1}(qt) X_1(0), \ldots, X_{m}(qt) -
\rho^{-1}(qt) X_{m}(0))
$$
we have by the triangle inequality
\BQNY
\abs{\X(qt)} \le
\abs{\X_{1/\rho}(qt)} + \frac1{\rho(qt)} \abs{\X(0)}.
\EQNY
Further, letting
$$
\zeta^{**}(qt) =\delta\abs{\X_{1/\rho}(qt)} +
\sqrt{1-\delta^2}\left(X_{m+1}(qt)-\rho^{-1}(qt)X_{m+1}(0)\right),
$$
and by utilising similar arguments as for $\zeta^*$ given in \eqref{def: zeta*}, we have that $\zeta^{**}(qt)$ is independent of $\zeta(qt)$ and $\zeta^{**}(qt)\equivdis \sqrt{1-\rho^2(qt)}/\rho(qt) \zeta(0)$. Therefore, for any $t$ \hH{positive}
satisfying $0 < t^{\alpha/2} < \lambda < \lambda_0$ and $u > 1$
\BQNY
\lefteqn{
 \pk{\zeta (qt) > u +\frac{\lambda}{u},   \zeta (0) \le u \Big\lvert \zeta(qt) > u }
 } \\
&& \le \pk{ \zeta^{**}(qt) > \frac{\lambda}{u} + u -  \frac{\zeta(0)}{\rho(qt)},   \zeta (0) \le u \Big\lvert  \zeta (qt) > u} \\
& & \le \pk{ \zeta^{**}(qt) > \frac{\lambda}{u}  - \left(\frac{1}{\rho(qt)} -1\right) u \Big\lvert  \zeta (qt) > u} \\
& &\le  \pk{ \zeta^{**}(qt) > \frac{\lambda}{2u}}  = \pk{ \zeta (0) >
\frac{\rho(qt)}{\sqrt{1-\rho^2(qt)}} \frac\lambda{2u}} \\
&& \le \pk{ \zeta (0) > \frac\lambda{8t^{\alpha/2}}}.
\EQNY
Consequently, \tn{by Chebyshev's inequality} for any positive
constant $p>2/\alpha$
\BQNY
\pk{ \zeta (qt) > u +\frac{\lambda}{u},  \zeta (0) \le u}
\le 8^p\E{\abs{ \zeta (0)}^p}t^{\alpha p/2}\lambda^{-p} \pk{ \zeta (qt)>u}
\EQNY
holds for any $t$ \hH{positive} satisfying $0< t^{\alpha/2} <
\lambda < \lambda_0$ and $u$ large. Thus the proof is complete.
\qed

\noindent\textbf{Proof of Theorem}~\ref{ThmA} With \nelem{LTA1}--\nelem{LTC}, we conclude that the claim follows by an application of \netheo{T1}. \qed

\subsection{Proof of Theorem \ref{ThmB}}\label{proof.th3}
In view of   \cite{Albin1990,Albin2003} or \cite{leadbetter1983extremes}, we need to verify two additional conditions  (\tn{see} Lemmas \ref{LTD} and \ref{LTD2}) for the order
statistics processes generated by the stationary Gaussian process $X$.
\begin{lemma} \label{LTD}
Under the assumptions of  \netheo{ThmB}, \hH{we have} \cling{for} any constants $a, T>0$
\BQN \label{cond: D}
 \lim_{\ve\downarrow 0}\limsup_{u\to\IF} \sum_{j =
[T/(aq)]}^{[\ve / \pk{X_{n:n}(0) >u}]} \pk{{X_{n:n}}(aqj)
>u \Big \lvert X_{n:n}(0) >u } = 0.
\EQN
\end{lemma}
\textbf{Proof}.
\tb{Recalling that  $X(t) - \rho(t)X(0)$ is independent of $X(0)$, }
 \BQN
 \lefteqn{
 \pk{X_{n:n}(t)>u \Big \lvert X_{n:n}(0)>u} = \pk{X_{n:n}(t)>u, X_{n:n}(0)>u \Big \lvert X_{n:n}(0)>u}
 } \notag\\
 &&= 2^n\left(  \pk{X(t)>X(0)>u \Big \lvert X(0)>u} \right)^n\notag\\
 &&\le  2^n\left(  \pk{X(t) - \rho(t)X(0) >u (1-\rho(t)), X(0) >u \Big \lvert X(0)>u} \right)^n \notag \\
 &&\le 2^n  \left(1- \Phi\left(u\sqrt{\frac{1- \abs{\rho(t)}}{1+\abs{\rho(t)}}}\right) \right)^n \notag\\
&& \le K u^{-n} \fracl{1- \cL{\abs{\rho(t)}}}{1+ \abs{\rho(t)}}^{-n/2} \expon{-\frac{nu^2}{2}\frac{1- \abs{\rho(t)}}{1+ \abs{\rho(t)}}}
\label{eqn: D}
 \EQN
holds for some positive constant $K$ and $u$ large (the constant $K$ below may be different from line to line), \cling{here $\Phi(\cdot)$ denotes the standard normal df}.

Now we choose a function  $ g=g(u)$ such that $\lim_{u\to\IF}g(u)=\IF,  {|\rho(g(u))|} = u^{-2}$. \hH{It} follows from  $u^{-2}\ln g(u) = o(1)$ that $g(u) \le \exp(\epsilon' u^2)$ for some $0<\epsilon'<\cling{n}/2(1-|\rho(T)|)/(1+|\rho(T)|)$ and sufficiently large $u$.
Now we split the sum in \eqref{cond: D} at $aqj = g(u)$. The first term satisfies
\BQNY
\lefteqn{
\sum_{j =[T/(aq)]}^{[g(u)/ (aq)]} \pk{X_{n:n}(aqj) >u \Big \lvert X_{n:n}(0) >u }}\\
&& \le  K\frac{g(u)}{aq} u^{-n} \fracl{1- \abs{\rho(T)}}{1+ \abs{\rho(T)}}^{-n/2} \expon{-\frac{nu^2}{2}\frac{1- \abs{\rho(T)}}{1+ \abs{\rho(T)}}}
 \\
& &\le {K  u^{2/\alpha -n}\expon{\epsilon' u^2-\frac{nu^2}{2}\frac{1- \abs{\rho(T)}}{1+ \abs{\rho(T)}}}}\\
&&\to 0, \quad u\to\IF
\EQNY
since $\epsilon'<n/2(1- \abs{\rho(T)})/(1+ \abs{\rho(T)})$. For the remaining term we have
\BQNY
\lefteqn{
\sum_{j = [g(u)/ (aq)]}^{[\ve /\pk{X_{n:n}(0) >u}]} \pk{X_{n:n}(aqj) >u \Big \lvert X_{n:n}(0) >u }}\\
&&\le  K \frac{\ve}{\pk{X_{n:n}(0) >u}} u^{-n} \fracl{1-u^{-2}}{1+u^{-2}}^{-n/2} \expon{-\frac{nu^2}{2}\frac{1-u^{-2}}{1+u^{-2}}} \\
&& \le  K \ve \expon{-\frac{nu^2}{2}\left(\frac{1-u^{-2}}{1+u^{-2}}-1\right)}\\
&& \le K \ve, \quad u\to\IF.
\EQNY
Therefore, the claim follows by taking $\ve\downarrow 0$. \qed

In the following lemma we shall establish the asymptotic independence of $X_{n:n}$
over suitable separate intervals (see condition $D'$ in \cite{Albin1990}).
\tn{In the notation used} below
$
\widetilde{\mathcal{A}}_{n,\alpha}$ is the constant appearing in \eqref{wideA}, and
\BQN
\label{def: T}
T= T(u) = \frac {(2 \pi)^{n/2}}{ \widetilde{\mathcal A}_{n,\alpha}} u^{n-\frac 2\alpha} \expon{ \frac{n u^2}2}.
\EQN

\begin{lemma} \label{LTD2}
Under the assumptions of \netheo{ThmB}, if futher $T=T(u)$ is defined by \eqref{def: T} and $a>0, 0 < \lambda < 1$ are given constants,
 then for any $0\le s_1 < \cdots < s_p < t_1 < \cdots < t_{p'}$ in
 $ \{ aqj: j\in \hH{\mathbb{Z}}, 0\le aqj \le T\} $ with $t_1 - s_p \ge \lambda T$ we   have
 \BQN \label{Asym.ind}
 \lefteqn{\hH{\limit{u}}\Bigg\lvert
 \pk{ \bigcap_{ i=1}^p\{{X_{n:n}}(s_i) \le u\}, \bigcap_{ j=1}^{p'}\{X_{n:n}(t_j) \le u\} } }\notag
 \\ &&\quad -
 \pk{\bigcap_{ i=1}^p\{X_{n:n}(s_i) \le u\}} \pk{\bigcap_{ j=1}^{p'}\{X_{n:n}(t_j) \le u\}}
 \Bigg\lvert    =0.
 \EQN
\end{lemma}
\textbf{Proof}. First, taking logarithms on both sides of \eqref{def: T} we
obtain
 \BQNY
&& \ln T = \frac{nu^2}2 + \left( n-\frac{2}{\alpha} \right)\ln u + \ln\hH{\Biggl(}\frac{(2\pi)^{n/2}}{\widetilde{\mathcal A}_{n,\alpha}}\Biggr),
 \EQNY
 which together with $u^2 = (2/n)\ln T (1+o(1))$ implies that
 \BQN  \label{Asym: T}
 u^2 = \frac{2\ln T}n + \left(\frac 2{n\alpha}-1\right)\ln \ln T +
 \ln \fracl n2^{1-\frac2{n\alpha}}{ \frac {(\widetilde{\mathcal{A}}_{n,\alpha})^{\frac2n}} {2\pi}} (1+o(1))
\EQN
as $T\to\IF$.
Further, define (hereafter  $\mathbb I\{\cdot\}$ denotes the indicator function)
\BQNY
X_{ij}= X_j(s_i)\mathbb I\{i\le p\} +  X_j(t_{i-p})\mathbb I\{p<i\le p+p'\}, \ \  1\le i\le p+p', 1\le j\le n,
\EQNY
and $\{Y_{ij}, 1\le i\le p, 1\le j\le n\}\equivdis \{X_{ij}, 1\le i\le p, 1\le j\le n\}$, independent of
$\{Y_{ij}, p+1\le i\le p+p', 1\le j\le n\}\equivdis \{X_{ij}, p+1\le i\le p+p', 1\le j\le n\}$.
Applying \nelem{MinP} with $X_{i(n)}= X_{n:n}(s_i)\mathbb I\{i\le p\} + X_{n:n}(t_{i-p})\mathbb I\{p<i\le p+p'\}$ \tb{(see the Appendix)}, using similar arguments as in Lemma 8.2.4 in \cite{leadbetter1983extremes}
we obtain that the left-hand side of \eqref{Asym.ind} is bounded from above by
\BQNY
\lefteqn{
K {u^{-2(n-1)}}\fracl Tq\sum_{\lambda T\le t_j-s_i \le T} e^{-\frac{{n u^2}}{1+\abs{\rho(t_j-s_i)}}}\int_0^{\abs{\rho(t_j-s_i)}}\frac{(1+{\abs h})^{2(n-1)}}{(1-h^2)^{n/2}}\, dh
} \\
&& \le Ku^{-2(n-1)}\fracl Tq \sum_{\lambda T \le aqj \le T} \abs{\rho(aqj)} e^{-\frac{nu^2}{1+\abs{\rho(aqj)} }}, \quad \mbox{for\ large\ } u,
\EQNY
where $K$ is some {positive} constant. \cL{The rest of the {proof} consists of the same arguments as that of
Lemma 12.3.1 in \cite{leadbetter1983extremes} by using \eqref{Asym: T} and the Berman\rq{}s condition $\rho(t)\ln t = o(1)$. Hence the proof is complete.}
\COM{
  (the constant $K$ below may \tb{be} different from line to line). \\
Next, letting $\delta(t) = \sup\{\abs{\rho(s)}\ln s: s\ge t\}, t\ge 1$,  we have that $\abs{\rho(t)} \le \delta(t) /\ln t$ and $\delta(t) \le M$ for some positive constant $M$ and all sufficiently large $t$.
Therefore, by \eqref{Asym: T}
\BQNY
\expon{{-\frac{nu^2}{1+\abs{\rho(aqj)}}}}
& \le&
\expon{-nu^2\left(1 - \frac{\delta(\lambda T)}{\ln(\lambda T)}\right)}
\\
&\le&  { K \expon{-nu^2}}\\
&\le& K T^{-2}(\ln T)^{n-2/\alpha}
\EQNY
holds for $T$ large.
Consequently,
\BQNY
\lefteqn{
 K  u^{-2(n-1)}\frac Tq\sum_{\lambda T \le aqj \le T} \abs{\rho(aqj)} \expon{{-\frac{nu^2}{1+\abs{\rho(aqj)}}}}
 } \\
 && \le  K u^{-2(n-1)}\fracl Tq^2\frac1{T/q}\sum_{\lambda T \le aqj \le T} \abs{\rho(aqj)} {\ln(aqj) }\frac1{\ln(\lambda T)}T^{-2}(\ln T)^{n-2/\alpha} \\
&& \le K \frac1{T/q}\sum_{\lambda T \le aqj \le T} \abs{\rho(aqj)} \ln(aqj), \quad T\to\IF
\EQNY
since {$\rho(t)\ln t = o(1)$.} Hence the proof is complete.}
\qed

\noindent\textbf{Proof of Theorem}~\ref{ThmB} Since \netheo{ThmA} and Lemmas \ref{LTD} and \ref{LTD2} hold for the $n$th order statistics process
$X_{n:n}$, in view of Lemma B in \cite{Albin2003} we have for $T=T(u)$ defined as in \eqref{def: T}
\BQNY
 \lim_{u\to\IF}\pk{\sup_{t\in[0, T(u)]}X_{n:n}(t) \le u + \frac x{nu} } = \expon{-e^{-x}}, \quad x\in\R.
\EQNY
\hH{Hence the proof follows by expressing $u$ in terms of $T$ \tb{as in} \eqref{Asym: T}.}
\qed

\section{Appendix}\label{sec5}
Let {$\vk{X}= (X_1 \ldot X_d)$ and $\vk{Y}=(Y_1 \ldot Y_d)$} be {two Gaussian random vectors with $N(0,1)$ components} and covariance matrices ${\Sigma^1}=(\sigma^{(1)}_{ij})$ and
${\Sigma^0}=(\sigma^{(0)}_{ij})$, respectively.  The most elaborated version of Berman's
inequality is due to \tn{Li and Shao} \cite{LiShao02}, where it is shown that
 for $\vk u= (u_1,\ldots,u_d) \in \R^d$ \cling{(hereafter the notation $\vk{x} \le \vk{y}$ for any $\vk x, \vk y \in\R^d$ means $x_i\le y_i$ for all $i\le d$)}
\BQN\label{ARIJ}
\Abs{\pk{\vk X\leq \vk u}-\pk{\vk Y\leq\vk u} }
\leq \frac{1}{2 \pi }\sum_{1\le i< j\le d}
\ARIJ
\expon{-\frac{u^2_{i}+u^2_{j}}{2(1+\rho_{ij})}},
\EQN
where
$ \rho_{ij}:=\max(|\sigma^{(1)}_{ij}|,|\sigma^{(0)}_{ij}|), \ARIJ:=|{\arcsin(\sigma^{(1)}_{ij})- \arcsin( \sigma^{(0)}_{ij})}|.$

Our goal is to establish Li and Shao's extension of comparison lemma (Berman's inequality) for the minimum and the maximum order statistics of Gaussian random vectors.
Let therefore $\vk{X}_j = ( X_{ij}, i\le d),  j\le n$ be $n$ independent copies of $\vk{X}$.
Denote the minimum and maximum order statistics vector $\X_{n:n} = (X_{i(n)}, i\le d), \X_{1:n} = (X_{i(1)}, i\le d)$ with $X_{i(n)}= \min_{j\le n} X_{ij},
X_{i(1)}= \max_{j\le n} X_{ij}, i\le d$.
\hH{Similarly, for $\vk{Y}_j, j\le n$ independent copies of $\vk{Y}$ we define
the minimum and maximum order statistics vectors $\vk Y_{n:n}$ and $\vk Y_{1:n}$, respectively}.

\begin{lemma}\label{Min}
\hH{For arbitrary $\vk u \in \R^d$ and \tbb{$k=1,n$} we have}
\BQN\label{eq:Min}
\ABs{
\pk{\X_{\tbb{k}:n}\le \vk u} -\pk{\vk Y_{\tbb{k}:n} \le \vk u}}  \le
\frac{n}{2\pi} \sum_{1\le i< j\le d} \ARIJ \expon{-\frac{u_i^2+u_j^2}{2(1+\rho_{ij})}}.\qquad
\EQN
\end{lemma}
\textbf{Proof}. Note that $- \vk{X}$ and $-\vk{Y}$ have the same distributions as those of $\vk{X}$ and $\vk{Y}$, respectively.
Using Theorem 2.1 in \cite{LuW2014} with constants $\lambda_{ij}=-u_i, i\le d, j\le n$, we have
\BQNY
\lefteqn{\Abs{ \pk{\X_{n:n} \le \vk u} -
\pk{\vk Y_{n:n} \le \vk u}}}\\
& &=
\Abs{\pk{\cup_{i=1}^d\cap_{j=1}^n\{-Y_{ij} \le -u_i\}} - \pk{\cup_{i=1}^d\cap_{j=1}^n\{-X_{ij} \le -u_i\}}}
\notag\\ &&\le  \frac n{2\pi}\sum_{1\le i<l\le d}A_{il}\expon{-\frac{u_i^2+u_l^2}{2(1+\rho_{il})}}.
\EQNY
Next, since $\abs{a^n-b^n} \le n\abs{a-b}, a,b\in[0,1]$ and $n\inn$, we have by \eqref{ARIJ}
\BQNY
\abs{\pk{\X_{1:n} \le \vk u} -
\pk{\vk Y_{1:n} \le \vk u}} &=&
\abs{\big(\pk{\vk X \hH{\le}\vk u}\big)^n - \big(\pk{\vk Y \hH{\le}\vk u}\big)^n
}\notag\\
&\le&  \frac n{2\pi}\sum_{1\le i<l\le d}A_{il} \expon{-\frac{u_i^2+u_l^2}{2(1+\rho_{il})}}, \qquad\qquad
\EQNY
hence \eqref{eq:Min} for \hH{$k=1$} follows and thus the proof is complete.  \qed

\begin{lemma}\label{MinP} Let $\X_{n:n}$ and $\vk Y_{n:n}$ be the minimum order statistics vectors defined above.
Then, {for all $\vk u >\vk 0$}
\BQNY
\ABs{
\pk{\X_{n:n}\le \vk u} -\pk{\vk Y_{n:n} \le \vk u}}  \le
{\frac n{(2\pi)^n} }u^{-2(n-1)}\sum_{1\le i< l\le d}\abs{A_{il}^*} \expon{-\frac{n u^2}{1+\rho_{il}}},
\EQNY
where $u=\min_{1\le i\le n}u_i$ and
\BQN\label{Ail*}
A_{il}^* = \int_{\sigma_{il}^{(0)}}^{\sigma_{il}^{(1)}}\frac{(1+{\abs h})^{2(n-1)}}{(1-h^2)^{n/2}}\, dh.
\EQN
\end{lemma}
\textbf{Proof}.
We follow the idea of the proof of {Theorem 2.1 in \cite{LuW2014}. Let $\{Z_{ij}^h, i\le d, j\le n\}$ be $N(0,1)$ \hH{rvs} with covariance matrix $\cL\Sigma^h =(\sigma^h_{ij,lk})$ where
$$
\sigma^h_{ij,lk} = \E{Z_{ij}^hZ_{lk}^h} = \tau^h_{il} \mathbb I\{j=k\},\quad i, l\le d,\ j,k\le n,\ h\in[0,1],
$$
with {$\tau^h_{il}:=h\sigma^{(1)}_{il} + (1-h)\sigma^{(0)}_{il}$}.
{Clearly,} $\vk Z_j^h= \{Z_{ij}^h, i\le d\}, j\le n$ are independent and identically $N_d(0, \Sigma^h)$ distributed with $\Sigma^h=h\Sigma^1 + (1-h)\Sigma^0$.
Without loss of generality, we assume that $\Sigma^1$ and $\Sigma^0$ are positive definite.
Consequently, we have (see (3.4) and (3.19) in \cite{LuW2014})
\BQNY
\lefteqn{\pk{\vk X_{n:n} \hH{\le}\vk u} -\pk{\vk Y_{n:n} \hH{\le}\vk u}}\notag \\
& = &\pk{{\vk Z_{1:n}^1} \ge - \vk u}- \pk{{\vk Z_{1:n}^0} \ge - \vk u} \notag \\
& = & n\sum_{1\le i<l\le d} (\sigma^{(1)}_{il} -\sigma^{(0)}_{il})\int_0^1 \, dh\Big(\varphi(-u_i, -u_l;  \tau ^h_{il})\notag \\
&&  \times\pk{\cap_{s=1, s\neq i,l}^d\{\hH{W^h_s} > -u_s\}\cap_{t=2}^n\{Z_{it}^h \le -u_i, Z_{lt}^h \le -u_l\} },
\EQNY
where $W^h_s = \max_{1\le j\le n}Z_{sj}^h$, and $\varphi(-u_i, -u_l\cling{;}   \tau^h_{il})$ is the bivariate pdf of $(Z_{i1}^h, Z_{l1}^h)$
which satisfies
\BQNY
\varphi(-u_i, -u_l;  \tau^h_{il}) &=& \frac{1}{2\pi (1-(\tau^h_{il})^2)^{1/2}} \expon{-\frac{u_i^2-2\tau^h_{il} u_iu_l +u_l^2 }{2(1-(\tau^h_{il})^2)}}
\\
&\le & \frac{1}{2\pi (1-(\tau^h_{il})^2)^{1/2}} \expon{-\frac{u^2}{1+\rho_{il}}}, \quad u=\min_{1\le i\le n}u_i.
\EQNY
Next, let \hH{$(Z_i, \tilde Z_l)$} be a bivariate standardized normal random vector with correlation $\lvert\tau_{il}^h\lvert$
 and set $u=\min_{1\le i\le n} u_i>0$.  \hH{Slepian's inequality in \cite{Pit96} and Lemma 2.3 in \cite{PicandsA} \cling{imply}}
\BQNY
\pk{Z_{it}^h\hH{\le}-u_i, Z_{lt}^h \hH{\le}-u_l}&\le& \pk{Z_i\hH{\le}-u_i, \cling{\tilde Z_l} \hH{\le}-u_l}\\
 & \le &\frac{(1+\lvert\tau_{il}^h\lvert)^2}{u^2}\varphi(u, u; \lvert\tau_{il}^h\lvert),
\quad t\le n.
\EQNY
Consequently, with $A_{il}^*$ defined by \eqref{Ail*} \cL{and $x_+=\max(x,0)$}
\BQNY
\lefteqn{\pk{\vk X_{n:n} \hH{\le}\vk u} -\pk{\vk Y_{n:n} \hH{\le}\vk u}}\notag \\
 && \le \frac n{(2\pi)^n}  u^{-2(n-1)}\sum_{1\le i<l\le d} (\sigma^{(1)}_{il} -\sigma^{(0)}_{il})_+
 \expon{-\frac{nu^2}{1+\rho_{il}}}\int_0^1 \frac{(1+\lvert\tau_{il}^h\lvert)^{2(n-1)} }{(1-(\tau^h_{il})^2)^{n/2}}\, dh
 \notag\\ && \le \frac n{(2\pi)^n}  u^{-2(n-1)}\sum_{1\le i<l\le d} (A_{il}^*)_+\expon{-\frac{nu^2}{1+\rho_{il}}},
\EQNY
 where in the last step we used the equality $\tau^h_{il}=h(\sigma_{il}^{(1)} - \sigma_{il}^{(0)})+ \sigma_{il}^{(0)}$, and thus
\BQNY
\int_0^1  \frac{(1+\lvert\tau _{il}^h\lvert)^{2(n-1)} }{(1-(\tau^h_{il})^2)^{n/2}}\, dh = \frac1{\sigma_{il}^{(1)} - \sigma_{il}^{(0)}}\int_{\sigma_{il}^{(0)}}^{\sigma_{il}^{(1)}}\frac{(1+\abs h)^{2(n-1)}}{(1-h^2)^{n/2}}\, dh.
 \EQNY
\tn{This completes the proof.}
\qed

\begin{acknowledgements}
\eH{We are grateful} to the referees for their careful reading
and numerous suggestions which greatly improved the paper.
\end{acknowledgements}

\def\polhk#1{\setbox0=\hbox{#1}{\ooalign{\hidewidth
  \lower1.5ex\hbox{`}\hidewidth\crcr\unhbox0}}}


\begin{thebibliography}{10}
\providecommand{\url}[1]{{#1}}
\providecommand{\urlprefix}{URL }
\expandafter\ifx\csname urlstyle\endcsname\relax
  \providecommand{\doi}[1]{DOI~\discretionary{}{}{}#1}\else
  \providecommand{\doi}{DOI~\discretionary{}{}{}\begingroup
  \urlstyle{rm}\Url}\fi

\bibitem{AlRawwashS07}
Al-Rawwash, M., Seif, M.: Measuring the surface roughness using the spatial
  statistics application.
\newblock J. Appl. Statist. Sci. \textbf{15}(2, [2006 on cover]), 205--213
  (2007)

\bibitem{AlbinPHD}
Albin, J.\eH{M.P.}: On extremal theory for non differentiable stationary processes.
\newblock PhD Thesis, University of Lund, Sweden  (1987)

\bibitem{Albin1990}
Albin, J.\eH{M.P.}: On extremal theory for stationary processes.
\newblock Ann. Probab. \textbf{18}(1), 92--128 (1990)

\bibitem{Albin1992}
Albin, J.\eH{M.P.}: Extremes and crossings for differentiable stationary processes with
  application to {G}aussian processes in {${\bf R}^m$} and {H}ilbert space.
\newblock Stochastic Process. Appl. \textbf{42}(1), 119--147 (1992)

\bibitem{Albin1998}
\cling{ Albin, J.M.P.: On extremal theory for self-similar processes.
\newblock  Ann. Probab. \textbf{26}(2), 743--793 (1998)}

\bibitem{Albin2003}
Albin, J.\eH{M.P.}, Jaru{\v{s}}kov{\'a}, D.: On a test statistic for linear trend.
\newblock Extremes \textbf{6}(3), 247--258 (2003)

\bibitem{MR2775212}
Alodat, M.: An approximation to cluster size distribution of two {G}aussian
  random fields conjunction with application to {FMRI} data.
\newblock J. Statist. Plann. Inference \textbf{141}(7), 2331--2347 (2011)

\bibitem{AlodatR09}
Alodat, M., Al-Rawwash, M.: Skew-{G}aussian random field.
\newblock J. Comput. Appl. Math. \textbf{232}(2), 496--504 (2009)

\bibitem{MR2654766}
Alodat, M., Al-Rawwash, M., Jebrini, M.: Duration distribution of the
  conjunction of two independent {$F$} processes.
\newblock J. Appl. Probab. \textbf{47}(1), 179--190 (2010)

\bibitem{Aue2009}
Aue, A., Horv{\'a}th, L., Hu{\v{s}}kov{\'a}, M.: Extreme value theory for
  stochastic integrals of {L}egendre polynomials.
\newblock J. Multivariate Anal. \textbf{100}(5), 1029--1043 (2009)

\bibitem{Berman82}
Berman, S.: Sojourns and extremes of stationary processes.
\newblock Ann. Probab. \textbf{10}(1), 1--46 (1982)

\bibitem{Berman92}
Berman, S.: Sojourns and extremes of stochastic processes.
\newblock Wadsworth \& Brooks/Cole Advanced Books \& Software, Pacific Grove,
  CA (1992)

\bibitem{Cambanis81}
Cambanis, S., Huang, S., Simons, G.: On the theory of elliptically contoured
  distributions.
\newblock J. Multivariate Anal. \textbf{11}(3), 368--385 (1981)

\bibitem{ChengXiao13}
Cheng, D., Xiao, Y.: Geometry and excursion probability of multivariate
  {G}aussian random fields.
\newblock \textbf{Manuscript} (2013)

\bibitem{DebickiHJminima}
D{\c{e}}bicki, K., Hashorva, E., Ji, L., Tabi{\'s}, K.: On the probability of
  conjunctions of stationary {G}aussian processes.
\newblock Statist. Probab. Lett. \textbf{88}, 141--148 (2014)

\bibitem{DiekerY}
\eH{Dieker, A.B., Yakir, B.: On asymptotic constants in the theory of {G}aussian processes.
\newblock Bernoulli \textbf{20}, 1600--1619 (2014)}


\bibitem{Faletal2010}
Falk, M., H\"usler, J., Reiss, R.D.: Laws of small numbers: {E}xtremes and rare
  events.
\newblock In: DMV Seminar, vol.~23, p. 3rd edn. Birkh\"auser, Basel (2010)

\bibitem{FarkasH13}
Farkas, Y., Hashorva, E.: Tail approximation for reinsurance portfolios of
  Gaussian-like risks.
\newblock Scand. Act. J. \textit{DO\cling{I} 10.1080/03461238.2013.825639} (2014)

\bibitem{FERET}
Freitas, A., H\"usler, J., Temido, M.: Limit laws for maxima of a stationary
  random sequence with random sample size.
\newblock TEST \textbf{21}(1), 116--131 (2012)

\bibitem{HJ14}
Hashorva, E., Ji, L.: Piterbarg theorems for chi-processes with trend.
\newblock Extremes \textit{DOI 10.1007/s10687-014-0201-1}  (2014)


\bibitem{leadbetter1983extremes}
Leadbetter, M., Lindgren, G., Rootz{\'e}n, H.: Extremes and related properties
  of random sequences and processes, vol.~11.
\newblock Springer Verlag (1983)

\bibitem{LiShao02}
Li, W., Shao, Q.M.: A normal comparison inequality and its applications.
\newblock Probab. Theory Related Fields \textbf{122}(4), 494--508 (2002)

\bibitem{LuW2014}
Lu, D., Wang, X.: Some new normal comparison inequalities related to {G}ordon's
  inequality.
\newblock Statist. Probab. Lett. \textbf{88}, \cL{133--140} (2014)

\bibitem{PicandsA}
Pickands III, J.: Upcrossing probabilities for stationary {G}aussian processes.
\newblock Trans. Amer. Math. Soc. \textbf{145}, 51--73 (1969)

\bibitem{Pit96}
Piterbarg, V.\eH{I}.: Asymptotic methods in the theory of {G}aussian processes and
  fields, \emph{Translations of Mathematical Monographs}, vol. 148.
\newblock American Mathematical Society, Providence, RI (1996)

\bibitem{Oleg05}
Seleznjev, O.: Asymptotic behavior of mean uniform norms for sequences of
  {G}aussian processes and fields.
\newblock Extremes \textbf{8}(3), 161--169 \cling{(2006) } 

\bibitem{TANSPA}
Tan, Z., Hashorva, E.: Exact asymptotics and limit theorems for supremum of
  stationary {$\chi$}-processes over a random interval.
\newblock Stochastic Process. Appl. \textbf{123}(8), 2983--2998 (2013)

\bibitem{TANWU}
Tan, Z., Wu, C.: Limit laws for the maxima of stationary chi-processes under
  random index.
\newblock TEST \textbf{in press}  (2014)

\bibitem{TurkmanA}
\hH{Turkman, K.F.:
Discrete and continuous time series extremes of stationary processes.
  \textit{Handbook of statistics Vol 30. Time Series Methods and Aplications. Eds. T.S. Rao, S.S. Rao and C.R. Rao. Elsevier},
 565--580 (2012)
}

\bibitem{TurkmanB}
\hH{ Turkman, K.F., Turkman, M.A.A., Pereira, J.M.: Asymptotic models and inference for extremes of
              spatio-temporal data. Extremes \textbf{13}(4), 375--397 (2010)
}


\bibitem{MR1747100}
Worsley, K., Friston, K.: A test for a conjunction.
\newblock Statist. Probab. Lett. \textbf{47}(2), 135--140 (2000)

\eH{\bibitem{Yakir}
Yakir, B.:  Extremes in Random Fields: A Theory and its Applications,
\emph{Higher Publication Press}, Wiley, New York (2013). }



\bibitem{ZareifardJ2013}
Zareifard, H., Jafari~Khaledi, M.: Non-{G}aussian modeling of spatial data
  using scale mixing of a unified skew {G}aussian process.
\newblock J. Multivariate Anal. \textbf{114}, 16--28 (2013)


\end{thebibliography}
\end{document}